\documentstyle[12pt]{article}

\catcode`\@=11
\long\def\@makefntext#1{ %\parindent 1em
\protect\noindent \hbox to 3.2pt {\hskip-.9pt
$^{{\ninerm\@thefnmark}}$\hfil}#1\hfill} %can be used

\def\thefootnote{\fnsymbol{footnote}}
 \def\@makefnmark{\hbox to 0pt{$^{\@thefnmark}$\hss}}  %original

\def\ps@myheadings{\let\@mkboth\@gobbletwo
\def\@oddhead{\hbox{} %\sl
\rightmark\hfil\ninerm\thepage}
\def\@oddfoot{}\def\@evenhead{\ninerm\thepage\hfil %\sl
\leftmark\hbox{}}\def\@evenfoot{}
\def\sectionmark##1{}\def\subsectionmark##1{}}

\def\CVector#1#2#3{
\left\vert\begin{array}{c}
                 #1 \\
                 #2
                \end{array}, #3\right>}

\def\AVector#1#2#3{\left|  #1, #2; #3 \right>}

\begin{document}

%----------------------------PROCSLA.STY---------------------------------------
\newcommand{\symbolfootnote}{\renewcommand{\thefootnote}
        {\fnsymbol{footnote}}}
\renewcommand{\thefootnote}{\fnsymbol{footnote}}
\newcommand{\alphfootnote}
        {\setcounter{footnote}{0}
         \renewcommand{\thefootnote}{\sevenrm\alph{footnote}}}

%------------------------------------------------------------------------------
%NEW DEFINED SECTION COMMANDS
\newcounter{sectionc}\newcounter{subsectionc}\newcounter{subsubsectionc}
\renewcommand{\section}[1] {\vspace{0.6cm}\addtocounter{sectionc}{1}
\setcounter{subsectionc}{0}\setcounter{subsubsectionc}{0}\noindent
        {\bf\thesectionc. #1}\par\vspace{0.4cm}}
\renewcommand{\subsection}[1] {\vspace{0.6cm}\addtocounter{subsectionc}{1}
        \setcounter{subsubsectionc}{0}\noindent
        {\it\thesectionc.\thesubsectionc. #1}\par\vspace{0.4cm}}
\renewcommand{\subsubsection}[1]
{\vspace{0.6cm}\addtocounter{subsubsectionc}{1}
        \noindent {\rm\thesectionc.\thesubsectionc.\thesubsubsectionc.
        #1}\par\vspace{0.4cm}}
\newcommand{\nonumsection}[1] {\vspace{0.6cm}\noindent{\bf #1}
        \par\vspace{0.4cm}}

%NEW MACRO TO HANDLE APPENDICES
\newcounter{appendixc}
\newcounter{subappendixc}[appendixc]
\newcounter{subsubappendixc}[subappendixc]
\renewcommand{\thesubappendixc}{\Alph{appendixc}.\arabic{subappendixc}}
\renewcommand{\thesubsubappendixc}
        {\Alph{appendixc}.\arabic{subappendixc}.\arabic{subsubappendixc}}

\renewcommand{\appendix}[1] {\vspace{0.6cm}
        \refstepcounter{appendixc}
        \setcounter{figure}{0}
        \setcounter{table}{0}
        \setcounter{equation}{0}
        \renewcommand{\thefigure}{\Alph{appendixc}.\arabic{figure}}
        \renewcommand{\thetable}{\Alph{appendixc}.\arabic{table}}
        \renewcommand{\theappendixc}{\Alph{appendixc}}
        \renewcommand{\theequation}{\Alph{appendixc}.\arabic{equation}}
%       \noindent{\bf Appendix \theappendixc. #1}\par\vspace{0.4cm}}
        \noindent{\bf Appendix \theappendixc #1}\par\vspace{0.4cm}}
\newcommand{\subappendix}[1] {\vspace{0.6cm}
        \refstepcounter{subappendixc}
        \noindent{\bf Appendix \thesubappendixc. #1}\par\vspace{0.4cm}}
\newcommand{\subsubappendix}[1] {\vspace{0.6cm}
        \refstepcounter{subsubappendixc}
        \noindent{\it Appendix \thesubsubappendixc. #1}
        \par\vspace{0.4cm}}

%------------------------------------------------------------------------------
%MARCO FOR ABSTRACT BLOCK
\def\abstracts#1{{
        \centering{\begin{minipage}{30pc}\tenrm\baselineskip=12pt\noindent
        \centerline{\tenrm ABSTRACT}\vspace{0.3cm}
        \parindent=0pt #1
        \end{minipage} }\par}}

%------------------------------------------------------------------------------
%NEW MACRO FOR BIBLIOGRAPHY
\newcommand{\bibit}{\it}
\newcommand{\bibbf}{\bf}
\renewenvironment{thebibliography}[1]
        {\begin{list}{\arabic{enumi}.}
        {\usecounter{enumi}\setlength{\parsep}{0pt}
%1.25cm IS STRICTLY FOR PROCSLA.TEX ONLY
\setlength{\leftmargin 1.25cm}{\rightmargin 0pt}
%0.52cm IS FOR NEW DATA FILES
%\setlength{\leftmargin 0.52cm}{\rightmargin 0pt}
         \setlength{\itemsep}{0pt} \settowidth
        {\labelwidth}{#1.}\sloppy}}{\end{list}}

%------------------------------------------------------------------------------
%FOLLOWING THREE COMMANDS ARE FOR 'LIST' COMMAND.
\topsep=0in\parsep=0in\itemsep=0in
\parindent=1.5pc

%LIST ENVIRONMENTS
\newcounter{itemlistc}
\newcounter{romanlistc}
\newcounter{alphlistc}
\newcounter{arabiclistc}
\newenvironment{itemlist}
        {\setcounter{itemlistc}{0}
         \begin{list}{$\bullet$}
        {\usecounter{itemlistc}
         \setlength{\parsep}{0pt}
         \setlength{\itemsep}{0pt}}}{\end{list}}

\newenvironment{romanlist}
        {\setcounter{romanlistc}{0}
         \begin{list}{$($\roman{romanlistc}$)$}
        {\usecounter{romanlistc}
         \setlength{\parsep}{0pt}
         \setlength{\itemsep}{0pt}}}{\end{list}}

\newenvironment{alphlist}
        {\setcounter{alphlistc}{0}
         \begin{list}{$($\alph{alphlistc}$)$}
        {\usecounter{alphlistc}
         \setlength{\parsep}{0pt}
         \setlength{\itemsep}{0pt}}}{\end{list}}

\newenvironment{arabiclist}
        {\setcounter{arabiclistc}{0}
         \begin{list}{\arabic{arabiclistc}}
        {\usecounter{arabiclistc}
         \setlength{\parsep}{0pt}
         \setlength{\itemsep}{0pt}}}{\end{list}}

%------------------------------------------------------------------------------
%FIGURE CAPTION
\newcommand{\fcaption}[1]{
        \refstepcounter{figure}
        \setbox\@tempboxa = \hbox{\tenrm Fig.~\thefigure. #1}
        \ifdim \wd\@tempboxa > 6in
           {\begin{center}
        \parbox{6in}{\tenrm\baselineskip=12pt Fig.~\thefigure. #1 }
            \end{center}}
        \else
             {\begin{center}
             {\tenrm Fig.~\thefigure. #1}
              \end{center}}
        \fi}

%TABLE CAPTION
\newcommand{\tcaption}[1]{
        \refstepcounter{table}
        \setbox\@tempboxa = \hbox{\tenrm Table~\thetable. #1}
        \ifdim \wd\@tempboxa > 6in
           {\begin{center}
        \parbox{6in}{\tenrm\baselineskip=12pt Table~\thetable. #1 }
            \end{center}}
        \else
             {\begin{center}
             {\tenrm Table~\thetable. #1}
              \end{center}}
        \fi}

%------------------------------------------------------------------------------
%ACKNOWLEDGEMENT: this portion is from John Hershberger
\def\@citex[#1]#2{\if@filesw\immediate\write\@auxout
        {\string\citation{#2}}\fi
\def\@citea{}\@cite{\@for\@citeb:=#2\do
        {\@citea\def\@citea{,}\@ifundefined
        {b@\@citeb}{{\bf ?}\@warning
        {Citation `\@citeb' on page \thepage \space undefined}}
        {\csname b@\@citeb\endcsname}}}{#1}}

\newif\if@cghi
\def\cite{\@cghitrue\@ifnextchar [{\@tempswatrue
        \@citex}{\@tempswafalse\@citex[]}}
\def\citelow{\@cghifalse\@ifnextchar [{\@tempswatrue
        \@citex}{\@tempswafalse\@citex[]}}
\def\@cite#1#2{{$\null^{#1}$\if@tempswa\typeout
        {IJCGA warning: optional citation argument
        ignored: `#2'} \fi}}
\newcommand{\citeup}{\cite}

%------------------------------------------------------------------------------
%FOR FNSYMBOL FOOTNOTE AND ALPH{FOOTNOTE}
\def\fnm#1{$^{\mbox{\scriptsize #1}}$}
\def\fnt#1#2{\footnotetext{\kern-.3em
        {$^{\mbox{\sevenrm #1}}$}{#2}}}

%------------------------------------------------------------------------------
\font\twelvebf=cmbx10 scaled\magstep 1
\font\twelverm=cmr10 scaled\magstep 1
\font\twelveit=cmti10 scaled\magstep 1
\font\elevenbfit=cmbxti10 scaled\magstephalf
\font\elevenbf=cmbx10 scaled\magstephalf
\font\elevenrm=cmr10 scaled\magstephalf
\font\elevenit=cmti10 scaled\magstephalf
\font\bfit=cmbxti10
\font\tenbf=cmbx10
\font\tenrm=cmr10
\font\tenit=cmti10
\font\ninebf=cmbx9
\font\ninerm=cmr9
\font\nineit=cmti9
\font\eightbf=cmbx8
\font\eightrm=cmr8
\font\eightit=cmti8

%----------------------START OF DATA FILE------------------------------

\centerline{\tenbf THE USE OF QUANTUM GROUPS IN NUCLEAR STRUCTURE PROBLEMS}
%\baselineskip=16pt
%\centerline{\tenbf MANUSCRIPT USING COMPUTER SOFTWARE}
\vspace{0.8cm}
\centerline{\tenrm Dennis BONATSOS}
\baselineskip=13pt
\centerline{\tenit ECT$^*$, Villa Tambosi, Strada delle Tabarelle 286}
\baselineskip=12pt
\centerline{\tenit I-38050 Villazzano (Trento), Italy}
\vspace{0.3cm}
\centerline{\tenrm C. DASKALOYANNIS}
\baselineskip=13pt
\centerline{\tenit Department of Physics, Aristotle University of Thessaloniki}
\baselineskip=12pt
\centerline{\tenit GR-54006 Thessaloniki, Greece}
\vspace{0.3cm}
%\centerline{\tenrm and}
%\vspace{0.3cm}
\centerline{\tenrm P. KOLOKOTRONIS, D. LENIS}
\baselineskip=13pt
\centerline{\tenit Institute of Nuclear Physics, NCSR ``Demokritos''}
\baselineskip=12pt
\centerline{\tenit GR-15310 Aghia Paraskevi, Attiki, Greece}
\vspace{0.9cm}
\abstracts{Various applications of quantum algebraic techniques in nuclear
structure physics, such as the su$_q$(2) rotator model and its extensions,
the use of deformed bosons in the description of pairing correlations, and
the construction of deformed exactly soluble models (Interacting Boson Model,
Moszkowski model) are briefly reviewed. Emphasis is put in the study of
the symmetries of the anisotropic quantum harmonic oscillator with rational
ratios of frequencies, which underly the structure of superdeformed and
hyperdeformed nuclei, the Bloch--Brink $\alpha$-cluster model and possibly
the shell structure in deformed atomic clusters.
}

\vfil
%\vspace{0.8cm}
\rm\baselineskip=14pt
\section{Introduction}

Quantum algebras $^{1,2}$ (also called quantum groups) are deformed versions
of the
usual Lie algebras, to which they reduce when the deformation parameter
$q$ is set equal to unity. Their use in physics became popular with the
introduction $^{3-5}$ of the $q$-deformed harmonic oscillator as a tool for
providing a boson realization of the quantum algebra su$_q$(2), although
similar mathematical structures had already been known $^{6,7}$.
Initially used for solving the quantum Yang--Baxter equation, quantum algebras
have subsequently found applications in several branches of physics, as, for
example, in the description of spin chains, squeezed states, rotational
and vibrational nuclear and molecular spectra, and in conformal
field theories. By now several kinds of generalized deformed oscillators
$^{8-12}$
and generalized deformed su(2) algebras $^{13-19}$  have been introduced.

Here we shall confine ourselves to applications of quantum algebras in nuclear
structure physics. A brief description will be given of the su$_q$(2) rotator
model $^{20-25}$ and its extensions $^{13,26}$,  of the use of deformed
oscillators in the description of pairing correlations $^{27-29}$,
and of the formulation of deformed exactly soluble
models (Interacting Boson Model $^{30}$, Moszkowski model $^{31-33}$).
The purpose of this
short review is to provide the  reader with references for further reading.
Subsequently, the symmetries of the anisotropic quantum harmonic oscillator
with rational ratios of frequencies will be considered in more detail,
since they are of current interest $^{34,35}$ in connection with
superdeformed and hyperdeformed nuclei $^{36,37}$,
$\alpha$-cluster configurations in light nuclei $^{38-40}$,
and possibly with deformed atomic clusters $^{41,42}$.

\section{ The su$_q$(2) rotator model}

 The first application of quantum algebras in nuclear physics was the use
of the deformed algebra su$_q$(2) for the description of the rotational
spectra of deformed $^{20,21}$ and superdeformed $^{22}$ nuclei.
The same technique has been
used for the description of rotational spectra of diatomic molecules $^{23}$.
The Hamiltonian of the $q$-deformed rotator is proportional to the
second order Casimir operator of the su$_q$(2) algebra. Its Taylor expansion
contains powers of $J(J+1)$ (where $J$ is the angular momentum), being
similar to the expansion provided by the Variable Moment of Inertia (VMI)
model. Furthermore, the deformation
parameter $\tau$ (with $q=e^{i\tau}$) has been found to correspond to
the softness parameter of the VMI model $^{21}$.

B(E2) transition probabilities have also been described in this framework
$^{24}$.
In this case the $q$-deformed Clebsch--Gordan coefficients are used instead
of the normal ones. (It should be noticed that the $q$-deformed angular
momentum theory has already been much developed $^{24}$.) The model predicts
an increase of the B(E2) values with angular momentum, while the rigid
rotator model predicts saturation. Some experimental results supporting
this prediction already exist $^{24}$.

\section{ Extensions of the su$_q$(2) model}

The su$_q$(2) model has been successful in describing rotational nuclear
spectra. For the description of vibrational and transitional nuclear
spectra it has been found $^{26}$ that $J(J+1)$ has to be replaced by $J(J+c)$.
The additional parameter $c$ allows for the description of nuclear
anharmonicities in a way similar to that of the Interacting Boson Model
(IBM) $^{43,44}$ and the Generalized Variable Moment of Inertia (GVMI) model
$^{45}$.

Another generalization is based on the use of the deformed algebra
su$_{\Phi}$(2) $^{13}$, which is characterized by a structure function $\Phi$.
The usual su(2) and su$_q$(2) algebras are obtained for specific choices
of the structure function $\Phi$. The su$_{\Phi}$(2) algebra has been
constructed so that its representation theory resembles as much as possible
the representation theory of the usual su(2) algebra. Using this technique
one can construct, for example, a rotator having the same spectrum as the
one given by the Holmberg--Lipas formula $^{46}$. In addition to the
generalized
deformed su(2) algebra, generalized deformed oscillators $^{8-12}$
have also been
introduced and found useful in many physical applications.

\section{ Pairing correlations}

It has been found $^{27}$ that correlated fermion pairs coupled to zero
angular
momentum in a single-$j$ shell behave approximately as suitably defined
$q$-deformed bosons. After performing the same boson mapping to a simple
pairing Hamiltonian, one sees that the pairing energies are also correctly
reproduced up to the same order. The deformation parameter used ($\tau
=\ln q$) is found to be inversely proportional to the size of the shell,
thus serving as a small parameter.

The above mentioned system of correlated fermion pairs can be described
{\sl exactly} by suitably defined generalized deformed bosons $^{28}$. Then
both the commutation relations are satisfied exactly and the pairing energies
are reproduced exactly. The spectrum of the appropriate generalized
deformed oscillator corresponds, up to first order perturbation theory,
to a harmonic oscillator with an $x^4$ perturbation.

\section{ $q$-deformed versions of nuclear models}

A $q$-deformed version of a two dimensional toy Interacting Boson Model
(IBM) has been developed $^{30}$, mainly for testing the ways in which
spectra and transition probabilities are influenced by the $q$-deformation. A
$q$-deformed version of the full IBM is under development, while a
$q$-deformed version of the vibron model, which uses the IBM techniques
in the case of molecules, has already been developed $^{47}$.

Furthermore a $q$-deformed version of the Moszkowski model has been developed
$^{31}$ and RPA modes have been studied $^{32}$ in it. A $q$-deformed
Moszkowski model
with cranking has also been studied $^{33}$ in the mean-field approximation.
It has been seen that the residual interaction simulated by the
$q$-deformation is felt more strongly by states with large $J_z$. The
possibility of using $q$-deformation in assimilating temperature effects is
under discussion.

\section{ Anisotropic quantum harmonic oscillator with rational ratios of
frequencies}

The symmetries of the 3-dimensional anisotropic quantum harmonic oscillator
with rational ratios of frequencies (RHO) are of high current interest in
nuclear physics, since they are
the basic symmetries $^{34,35}$ underlying the structure of superdeformed and
hyperdeformed nuclei $^{36,37}$.
The 2-dimensional RHO is also of interest, in
connection with ``pancake'' nuclei, i.e. very oblate nuclei $^{35}$. Cluster
configurations in light nuclei can also be  described in terms of RHO
symmetries $^{38,39}$, which underlie the geometrical structure of the
Bloch--Brink $\alpha$-cluster model $^{40}$. The 3-dim RHO is also of
interest for the interpretation
of the observed shell structure in atomic clusters $^{41}$, especially
after the
realization $^{42}$ that large deformations can occur in such systems.

The two-dimensional $^{48-53}$ and
three-dimensional $^{54-60}$  anisotropic harmonic
oscillators have been the subject of several investigations, both at the
classical and the quantum mechanical level. These oscillators are examples
of superintegrable systems $^{61}$. The special cases with frequency
ratios 1:2 $^{62,63}$ and 1:3 $^{64}$ have also been considered. While
at the classical level it is clear that the su(N) or sp(2N,R) algebras can
be used for the description of the N-dimensional anisotropic oscillator, the
situation at the quantum level, even in the two-dimensional case, is not as
simple.

In this section  we are going to prove that a generalized deformed u(2)
algebra is the symmetry algebra of the two-dimensional anisotropic quantum
harmonic oscillator, which is the oscillator describing the single-particle
level spectrum of ``pancake'' nuclei, i.e. of triaxially deformed nuclei
with $\omega_x >> \omega_y$, $\omega_z$ $^{35}$.

\subsection{ The deformed u(2) algebra}

Let us consider the system described by the Hamiltonian:
\begin{equation}
\label{eq:Hamiltonian}H=\frac{1}{2}\left( {p_x}^2 + {p_y}^2 + \frac{x^2}{m^2}
+ \frac{y^2}{n^2} \right),
\end{equation}
where $m$ and $n$ are two natural numbers mutually prime ones, i.e. their
great common divisor is $\gcd (m,n)=1$.

We define the creation and annihilation operators $^{48}$
\begin{equation}
\label{eq:operators}
\begin{array}{ll}
a^\dagger=\frac{x/m - i p_x}{\sqrt{2}}, & a =
\frac{x/m + i p_x}{\sqrt{2}}, \\ [0.24in] b^\dagger=\frac{y/n - i p_y}{\sqrt{%
2}}, & b=\frac{y/n + i p_y}{\sqrt{2}}.
\end{array}
\end{equation}
These operators satisfy the commutation relations:
\begin{equation}
\label{eq:commutators}\left[ a,a^\dagger \right] = \frac{1}{m}, \quad \left[
b,b^\dagger \right] = \frac{1}{n}, \quad \mbox{other commutators}=0.
\end{equation}
Using Eqs (\ref{eq:operators}) and (\ref{eq:commutators}) we can prove by
induction that:
$$
\begin{array}{ll}
\left[ a, \left( a^\dagger \right)^p \right] = \frac{p}{m} \left( a^\dagger
\right)^{p-1} , & \left[ b, \left( b^\dagger \right)^p \right] =
\frac{p}{n} \left( b^\dagger \right)^{p-1} , \\ [0.24in] \left[ a^\dagger,
\left( a \right)^p \right] = - \frac{p}{m} \left( a \right)^{p-1} , & \left[
b^\dagger, \left( b \right)^p \right] = - \frac{p}{n} \left( b \right)^{p-1}
{}.
\end{array}
$$
Defining
$$
U=\frac{1}{2} \left\{ a, a^\dagger \right\}, \qquad W=\frac{1}{2} \left\{ b,
b^\dagger \right\},%
$$
one can easily prove that:
$$
\begin{array}{ll}
\left[ U, \left(a^\dagger \right)^p \right]= \frac{p}{m} \left(a^\dagger
\right)^p, & \left[ W, \left(b^\dagger \right)^p \right]=
\frac{p}{n} \left(b^\dagger \right)^p, \\ [0.24in] \left[ U, \left(a
\right)^p \right]= - \frac{p}{m} \left(a \right)^p, & \left[ W, \left(b
\right)^p \right]= - \frac{p}{n} \left(b \right)^p.
\end{array}
$$
Using the above properties we can define the enveloping algebra generated by
the operators:
\begin{equation}
\label{eq:generators}
\begin{array}{c}
S_+= \left(a^\dagger\right)^m \left(b\right)^n,\quad S_-= \left(a\right)^m
\left(b^\dagger\right)^n, \\
[0.24in] S_0= \frac{1}{2}\left( U - W \right), \quad H=U+W.
\end{array}
\end{equation}
These genarators satisfy the following relations:
\begin{equation}
\label{eq:SS}\left[ S_0,S_\pm \right]=\pm S_\pm, \quad \left[H,S_i\right]=0,
\quad \mbox{for}\quad i=0,\pm,
\end{equation}
and
$$
S_+S_- = \prod\limits_{k=1}^{m}\left( U - \frac{2k-1}{2m} \right)
\prod\limits_{\ell=1}^{n}\left( W + \frac{2\ell-1}{2n} \right),
$$
$$
S_-S_+ = \prod\limits_{k=1}^{m}\left( U + \frac{2k-1}{2m} \right)
\prod\limits_{\ell=1}^{n}\left( W - \frac{2\ell-1}{2n} \right).
$$
The fact that the operators $S_i$, $i=0, \pm$ are integrals of motion has
been already realized in $^{48}$.

The above relations mean that the harmonic oscillator of Eq. (\ref
{eq:Hamiltonian}) is described by the enveloping algebra of the
generalization of the u(2) algebra formed by the generators $S_0$, $S_+$, $%
S_-$ and $H$, satisfying the commutation relations of Eq. (\ref{eq:SS}) and
\begin{equation}
\label{eq:U2}
\begin{array}{c}
\left[S_-,S_+\right] = F_{m,n} (H,S_0+1) - F_{m,n} (H,S_0), \\
[0.24 in] \mbox{where}\quad F_{m,n}(H,S_0)= \prod\limits_{k=1}^{m}\left(
H/2+S_0 - \frac{2k-1}{2m} \right) \prod\limits_{\ell=1}^{n}\left( H/2-S_0 +
\frac{2\ell-1}{2n} \right).
\end{array}
\end{equation}
In the case of $m=1$, $n=1$ this algebra is the usual u(2) algebra, and the
operators $S_0,S_\pm$ satisfy the commutation relations of the ordinary u(2)
algebra, since in this case one easily finds that
$$
[S_-, S_+]=-2 S_0.%
$$
In the rest of the cases, the algebra is a deformed version of u(2), in
which the commutator $[S_-,S_+]$ is a polynomial of $S_0$ of order $m+n-1$.
In the case with $m=1$, $n=2$ one has
$$
[S_-,S_+]= 3 S_0^2 - H S_0 -{\frac{H^2}{4}} +{\frac{3}{16}},%
$$
i.e. a polynomial quadratic in $S_0$ occurs, while in the case of $m=1$, $n=3
$ one finds
$$
[S_-, S_+]= -4 S_0^3 + 3 H S_0^2 -{\frac{7}{9}} S_0 -{\frac{H^3}{4}} + {%
\frac{H}{4}},%
$$
i.e. a polynomial cubic in $S_0$ is obtained.

\subsection{ The representations}

The finite dimensional representation modules  of this algebra can be found
using the concept of the generalized deformed oscillator $^8$, in a
method similar to the one used in Ref. 65  for the study of quantum
superintegrable systems. The operators:
\begin{equation}
\label{eq:alge-gen}{\cal A}^\dagger= S_+, \quad {\cal A}= S_-, \quad {\cal N}%
= S_0-u, \quad u=\mbox{ constant},
\end{equation}
where $u$ is a constant to be determined, are the generators of a deformed
oscillator algebra:
$$
\left[ {\cal N} , {\cal A}^\dagger \right] = {\cal A}^\dagger, \quad \left[
{\cal N} , {\cal A} \right] = -{\cal A}, \quad {\cal A}^\dagger{\cal A}
=\Phi( H, {\cal N} ), \quad {\cal A}{\cal A}^\dagger =\Phi( H, {\cal N}+1 ).
$$
The structure function $\Phi$ of this algebra is determined by the function $%
F_{m,n}$ in Eq. (\ref{eq:U2}):
\begin{equation}
\label{eq:sf}
\begin{array}{l}
\Phi( H,
{\cal N} )= F_{m,n} (H,{\cal N} +u ) = \\ = \prod\limits_{k=1}^{m}\left( H/2+%
{\cal N} +u - \frac{2k-1}{2m} \right) \prod\limits_{\ell=1}^{n}\left( H/2-%
{\cal N} - u + \frac{2\ell-1}{2n} \right).
\end{array}
\end{equation}
The deformed oscillator corresponding to the structure function of Eq.  (\ref
{eq:sf}) has an energy dependent Fock space of dimension $N+1$ if
\begin{equation}
\label{eq:equations}\Phi(E,0)=0, \quad \Phi(E, N+1)=0, \quad \Phi(E,k)>0,
\quad \mbox{for} \quad k=1,2,\ldots,N.
\end{equation}
The Fock space is defined by:
\begin{equation}
H\vert E, k > =E \vert E, k >, \quad {\cal N} \vert E, k >= k \vert E, k
>,\quad a\vert E, 0 >=0,
\end{equation}
\begin{equation}
{\cal A}^\dagger \vert E, k> = \sqrt{\Phi(E,k+1)} \vert E, k+1>, \quad {\cal %
A} \vert E, k> = \sqrt{\Phi(E,k)} \vert E, k-1>.
\end{equation}
The basis of the Fock space is given by:
$$
\vert E, k >= \frac{1}{\sqrt{[k]!}} \left({\cal A}^\dagger\right)^k\vert E,
0 >, \quad k=0,1,\ldots N,
$$
where the ``factorial'' $[k]!$ is defined by the recurrence relation:
$$
[0]!=1, \quad [k]!=\Phi(E,k)[k-1]! \quad .
$$
Using the Fock basis we can find the matrix representation of the deformed
oscillator and then the matrix representation of the algebra of Eqs (\ref
{eq:SS}), (\ref{eq:U2}). The solution of Eqs (\ref{eq:equations}) implies
the following pairs of permitted values for the energy eigenvalue $E$ and
the constant $u$:
\begin{equation}
\label{eq:E1}E=N+\frac{2p-1}{2m}+\frac{2q-1}{2n} ,
\end{equation}
where $p=1,2,\ldots,m$, $q=1,2,\ldots,n$, and
$$
u=\frac{1}{2}\left( \frac{2p-1}{2m}-\frac{2q-1}{2n} -N \right),
$$
the corresponding structure function being given by:
\begin{equation}
\label{eq:structure-function}
\begin{array}{l}
\Phi(E,x)=\Phi^{N}_{(p,q)}(x)= \\
=\prod\limits_{k=1}^{m}\left( x +
\displaystyle \frac{2p-1}{2m}- \frac{2k-1}{2m} \right)
\prod\limits_{\ell=1}^{n}\left( N-x+ \displaystyle \frac{2q-1}{2n} + \frac{%
2\ell-1}{2n}\right) \\ =\displaystyle\frac{1}{m^m n^n} \displaystyle\frac{
\Gamma\left(mx+p\right) }{\Gamma\left(mx+p-m\right)} \displaystyle \frac{%
\Gamma\left( (N-x)n + q + n \right)} {\Gamma\left( (N-x)n + q \right)}.
\end{array}
\end{equation}
In all these equations one has $N=0,1,2,\ldots$, while the dimensionality of
the representation is given by $N+1$. Eq. (\ref{eq:E1})  means that there
are $m\cdot n$ energy eigenvalues corresponding to each $N$ value, each
eigenvalue having degeneracy $N+1$. (Later we shall see that the degenerate
states corresponding to the same eigenvalue can be labelled by an ``angular
momentum''.)

It is useful to show at this point that a few special cases are in agreement
with results already existing in the literature.

i) In the case $m=1$, $n=1$ Eq. (13) gives
$$ \Phi(E,x)= x(N+1-x),$$
while Eq. (12) gives
$$ E=N+1,$$
in agreement with Sec. IV.A of Ref. 65.

ii) In the case $m=1$, $n=2$ one obtains for $q=2$
$$ \Phi(E,x)= x(N+1-x)\left(N+{3\over 2}-x\right), \qquad E=N+{5\over 4},$$
while for $q=1$ one has
$$ \Phi(E,x)= x(N+1-x)\left(N+{1\over 2}-x\right), \qquad E=N+{3\over 4}.$$
These are in agreement with the results obtained in Sec. IV.F of Ref. 65
for the Holt potential (for $\delta =0$).

iii) In the case $m=1$, $n=3$ one has for $q=1$
$$\Phi(E,x)=x(N+1-x)\left(N+{1\over 3}-x\right) \left(N+{2\over 3}-x\right),
\qquad E=N+{2\over 3},$$
while for $q=2$ one obtains
$$\Phi(E,x)=x(N+1-x)\left(N+{2\over 3}-x\right) \left(N+{4\over 3}-x\right),
\qquad E=N+1,$$
and for $q=3$ one gets
$$\Phi(E,x)=x(N+1-x)\left(N+{4\over 3}-x\right) \left(N+{5\over 3}-x\right),
\qquad E=N+{4\over 3}. $$
These are in agreement with the results obtained in Sec. IV.D of Ref. 65
for the Fokas--Lagerstrom potential.

In all of the above cases we remark that the structure function has
the form
$$\Phi(x)= x(N+1-x) (\lambda+\mu x +\nu x^2 +\rho x^3+\sigma x^4+\dots,$$
which corresponds to a generalized deformed parafermionic algebra $^{66}$
 of order $N$, if $\lambda$, $\mu$, $\nu$, $\rho$, $\sigma$, \dots,
are real constants satisfying the conditions
$$\lambda + \mu x + \nu x^2 +\rho x^3 + \sigma x^4 +\dots > 0, \qquad
x\in \{1, 2, \dots, N\}.$$
These conditions are indeed satisfied in all cases.
It is easy to see that the obtained algebra corresponds to this of the
generalized parafermionic oscillator in all cases with frequency
ratios $1:n$.

The energy formula can be corroborated by using the
corresponding Schr\"{o}dinger equation. For the Hamiltonian of Eq. (\ref
{eq:Hamiltonian}) the eigenvalues of the Schr\"{o}dinger equation are given
by:
\begin{equation}
\label{eq:E2}E=\frac{1}{m}\left(n_x+\frac{1}{2}\right)+ \frac{1}{n}\left(n_y+%
\frac{1}{2}\right),
\end{equation}
where $n_x=0,1,\ldots$ and $n_y=0,1,\ldots$. Comparing Eqs (\ref{eq:E1}) and
(\ref{eq:E2}) one concludes that:
$$
N= \left[n_x/m\right]+\left[n_y/n\right],%
$$
where $[x]$ is the integer part of the number $x$, and
$$
p=\mbox{mod}(n_x,m)+1, \quad q=\mbox{mod}(n_y,n)+1.
$$

The eigenvectors of the Hamiltonian can be parametrized by the
dimensionality of the representation $N$, the numbers $p,q$, and the number $%
k=0,1,\ldots,N$. $k$ can be identified as $[n_x/m]$. One then has:

\begin{equation}
\label{eq:en-rep}H\left\vert
\begin{array}{c}
N \\
(p,q)
\end{array}
, k \right>= \left(N+\displaystyle
\frac{2p-1}{2m}+\frac{2q-1}{2n} \right)\left\vert
\begin{array}{c}
N \\
(p,q)
\end{array}
, k \right>,
\end{equation}
\begin{equation}
\label{eq:s0-rep}S_0 \left\vert
\begin{array}{c}
N \\
(p,q)
\end{array}
, k \right>= \left( k+ \displaystyle
\frac{1}{2} \left( \frac{2p-1}{2m}- \frac{2q-1}{2n} -N \right) \right)
\left\vert
\begin{array}{c}
N \\
(p,q)
\end{array}
, k \right>,
\end{equation}
\begin{equation}
\label{eq:sp-rep}S_+\left\vert
\begin{array}{c}
N \\
(p,q)
\end{array}
, k \right> = \sqrt{ \Phi^N_{(p,q)}(k+1)} \left\vert
\begin{array}{c}
N \\
(p,q)
\end{array}
, k +1\right>,
\end{equation}
\begin{equation}
\label{eq:sm-rep}S_-\left\vert
\begin{array}{c}
N \\
(p,q)
\end{array}
, k \right> = \sqrt{ \Phi^N_{(p,q)}(k)} \left\vert
\begin{array}{c}
N \\
(p,q)
\end{array}
, k -1\right>.
\end{equation}

\subsection{ The ``angular momentum'' quantum number}

It is worth noticing that the operators $S_0,S_\pm$ do not correspond to a
generalization of the angular momentum, $S_0$ being the operator
corresponding to the Fradkin operator $S_{xx}-S_{yy}$ $^{67,68}$.
The corresponding ``angular momentum'' is defined by:
\begin{equation}
\label{eq:angular-momentum}L_0=-i\left(S_+-S_-\right).
\end{equation}
The ``angular momentum'' operator commutes with the Hamiltonian:
$$
\left[ H,L_0 \right]=0.
$$
Let $\vert \ell> $ be the eigenvector of the operator $L_0$ corresponding to
the eigenvalue $\ell$. The general form of this eigenvector can be given by:
\begin{equation}
\vert \ell > = \sum\limits_{k=0}^N \frac{i^k c_k}{\sqrt{[k]!}} \left\vert
\begin{array}{c}
N \\
(p,q)
\end{array}
, k \right>.
\label{eq:ang-vector}
\end{equation}

In order to find the eigenvalues of $L$ and the coefficients $c_k$ we use
the Lanczos algorithm $^{69}$, as formulated in $^{70}$. From Eqs (\ref
{eq:sp-rep}) and (\ref{eq:sm-rep}) we find
$$
\begin{array}{l}
L_0|\ell >=\ell |\ell >=\ell \sum\limits_{k=0}^N\frac{i^kc_k}{\sqrt{[k]!}}%
\left|
\begin{array}{c}
N \\
(p,q)
\end{array}
,k\right\rangle = \\
=\frac 1i\sum\limits_{k=0}^{N-1}\frac{i^kc_k\sqrt{\Phi _{(p,q)}^N(k+1)}}{%
\sqrt{[k]!}}\left|
\begin{array}{c}
N \\
(p,q)
\end{array}
,k+1\right\rangle -\frac 1i\sum\limits_{k=1}^N\frac{i^kc_k\sqrt{\Phi
_{(p,q)}^N(k)}}{\sqrt{[k]!}}\left|
\begin{array}{c}
N \\
(p,q)
\end{array}
,k-1\right\rangle
\end{array}
$$
{}From this equation we find that:
$$
c_k=(-1)^k 2^{-k/2}H_k(\ell /\sqrt{2})/{\cal N},
\quad {\cal N}^2= \sum\limits_{n=0}^N 2^{-n}H_n^2(\ell /\sqrt{2})
$$
where the function $H_k(x)$ is a generalization of the ``Hermite''
polynomials (see also $^{71,72}$), satisfying the recurrence
relations:
$$
H_{-1}(x)=0,\quad H_0(x)=1,
$$
$$
H_{k+1}(x)=2xH_k(x)-2\Phi _{(p,q)}^N(k)H_{k-1}(x),
$$
and the ``angular momentum'' eigenvalues $\ell $ are the roots of the
polynomial equation:
\begin{equation}
H_{N+1}(\ell /\sqrt{2})=0.
\label{eq:eigenvalues}
\end{equation}
Therefore for a given value of $N$ there are $N+1$ ``angular momentum''
eigenvalues $\ell $, symmetric around zero (i.e. if $\ell $ is an ``angular
momentum'' eigenvalue, then $-\ell $ is also an ``angular momentum''
eigenvalue). In the case of the symmetric harmonic oscillator ($m/n=1/1$)
these eigenvalues are uniformly distributed and differ by 2. In the general
case the ``angular momentum'' eigenvalues are non-uniformly distributed. For
small values of $N$ analytical formulae for the ``angular momentum''
eigenvalues can be found $^{71}$. Remember that to each value of $N$
correspond $m\cdot n$ energy levels, each with degeneracy $N+1$.

In order to have a formalism corresponding to the one of the isotropic
oscillator, let us introduce  for every $N$ and
$(p,q)$ an ordering of the ``angular momentum'' eigenvalues
$$
\ell_m^{L,(p,q)}, \quad \mbox{where} \quad L=N
\quad \mbox{and} \quad m=-L,-L+2,\ldots,L-2,L,
$$
by assuming that:
$$
\ell_m^{L,(p,q)} \le \ell_n^{L,(p,q)} \quad \mbox{if} \quad m<n,
$$
 the corresponding eigenstate being given by:
\begin{equation}
\AVector{L}{m}{(p,q)}=
\sum\limits_{k=0}^N
\frac{(-i)^kH_k(\ell_m^{L,(p,q)} /\sqrt{2})}
{{\cal N} \sqrt{2^{k/2}[k]!}}
\CVector{N}{(p,q)}{k}.
\label{eq:harmonic}
\end{equation}
The above vector elements constitute  the analogue corresponding
to the  basis of ``sphe\-rical harmonic'' functions of the usual oscillator.

\subsection{ Examples}

For illustrative purposes, let us discuss a couple of examples in more
detail.

i) In the case $m=1$, $n=1$ the only $(p,q)$ value allowed is $(1,1)$. In
the Cartesian notation $|n_x n_y>$ the lowest energy corresponds to the
state $|00>$, the next energy level corresponds to the two states $|10>$ and
$|01>$,
while the next energy corresponds to the three states $|20>$, $|11>$, $|02>$.
 In the notation used in Eqs (15)-(18) these states can be written
as
$$ N(p,q)= 0(1,1), \quad E=1 \quad
  \rightarrow  \quad |00> =\CVector{0}{(1,1)}{0},  $$
$$
N(p,q)= 1(1,1), \quad E= 2 \quad  \rightarrow \quad
\left\{
\begin{array}{c}
|01>=\CVector{1}{(1,1)}{0}\\
|10>=\CVector{1}{(1,1)}{1},
\end{array}
\right.
$$
$$
N(p,q)= 2(1,1), \quad E=3 \quad  \rightarrow \quad
\left\{
\begin{array}{c}
|02>=\CVector{2}{(1,1)}{0}\\
|11>=\CVector{2}{(1,1)}{1} \\
|20>=\CVector{2}{(1,1)}{2}.
\end{array}
\right.
$$
{}From this example it is clear that the irreducible representations (irreps)
are characterized by the quantum numbers $N$ and $(p,q)$, while $k$
enumerates the degenerate states within each irrep.  The lowest irrep,
characterized by $N=0$ and $(p,q)=(1,1)$, has dimension $d=1$, while the
next irrep has $N=1$, $(p,q)=(1,1)$ and $d=2$ and the third one $N=2$,
$(p,q)=(1,1)$ and $d=3$.

Using the ``angular momentum'' eigenvalues defined by Eqs
(\ref{eq:eigenvalues}) and (\ref{eq:harmonic}) we find:
$$
E=1 \quad \rightarrow \quad
\AVector{0}{0}{(1,1)}=|0,0>, \quad  \leftarrow \quad \ell_0^{0,(1,1)}=0,
$$
$$
E=2 \quad \rightarrow \quad
\left\{
\begin{array}{l l}
\AVector{1}{-1}{(1,1)}=\frac{1}{\sqrt{2}}\left(|01>+i|1,0>\right)
&  \leftarrow \quad \ell_{(-1)}^{1,(1,1)}=-1\\
\AVector{1}{+1}{(1,1)}=\frac{1}{\sqrt{2}}\left(|01>-i|1,0>\right)
& \leftarrow \quad \ell_{(+1)}^{1,(1,1)}=1,
\end{array}
\right.
$$
$$
E=3 \quad \rightarrow \quad
\left\{
\begin{array}{l l}
\AVector{2}{-2}{(1,1)}=
\frac{1}{2}|0,2>+\frac{i}{\sqrt{2}}|1,1>-\frac{1}{2}|2,0>
&  \leftarrow \quad \ell_{(-2)}^{2,(1,1)}=-2\\
\AVector{2}{0}{(1,1)}=
\frac{1}{\sqrt{2}}|0,2>+\frac{1}{\sqrt{2}}|2,0>
&  \leftarrow \quad \ell_{(0)}^{2,(1,1)}=0\\
\AVector{2}{+2}{(1,1)}=
\frac{1}{2}|0,2>-\frac{i}{\sqrt{2}}|1,1>-\frac{1}{2}|2,0>
&  \leftarrow \quad \ell_{(+2)}^{2,(1,1)}=+2 .
\end{array}
\right.
$$
We see therefore that the angular momentum eigenvalues $l$ can be used,
instead of $k$, for enumerating the degenerate states within each irrep.

ii) In the case $m=1$, $n=2$ the allowed values of $(p,q)$ are $(1,1)$ and $%
(1,2)$. The lowest irrep is characterized by $N=0$, $(p,q)=(1,1)$, has
dimension $d=1$ and contains the state $|n_x,n_y>=|00>$, the next irrep is
characterized by $N=0$, $(p,q)=(1,2)$, has $d=1$ and contains the state $|01>
$, the third irrep has $N=1$, $(p,q)=(1,1)$, $d=2$ and contains the states $
|10>$ and $|02>$, the fourth irrep has $N=1$, $(p,q)=(1,2)$, $d=2$ and
contains the states $|11>$ and $|03>$, the fifth irrep has $N=2$, $(p,q)=
(1,1)$, $d=3$ and contains the states $|20>$, $|12>$ and $|04>$, the sixth
irrep has $N=2$, $(p,q)=(1,2)$, $d=3$ and contains the states $|21>$, $|13>$
and $|05>$. The states are listed in both notations below:

$$
 E=\frac{3}{4}
 \quad  \rightarrow  \quad |00> =\CVector{0}{(1,1)}{0},
\qquad E=\frac{5}{4}
 \quad  \rightarrow  \quad |01> =\CVector{0}{(1,2)}{0},
$$
$$
E= \frac{7}{4} \quad  \rightarrow \quad
\left\{
\begin{array}{c}
|02>=\CVector{1}{(1,1)}{0}\\
|10>=\CVector{1}{(1,1)}{1},
\end{array}
\right.
\qquad E= \frac{9}{4} \quad  \rightarrow \quad
\left\{
\begin{array}{c}
|03>=\CVector{1}{(1,2)}{0}\\
|11>=\CVector{1}{(1,2)}{1},
\end{array}
\right.
$$
$$
E=\frac{11}{4}  \rightarrow
\left\{
\begin{array}{c}
|04>=\CVector{2}{(1,1)}{0}\\
|12>=\CVector{2}{(1,1)}{1} \\
|20>=\CVector{2}{(1,1)}{2},
\end{array}
\right. \qquad
E=\frac{13}{4}  \rightarrow
\left\{
\begin{array}{c}
|05>=\CVector{2}{(1,2)}{0}\\
|13>=\CVector{2}{(1,2)}{1} \\
|21>=\CVector{2}{(1,2)}{2}.
\end{array}
\right.
$$

Using the ``angular momentum'' eigenvalues defined by Eqs
(\ref{eq:eigenvalues}) and (\ref{eq:harmonic}) we find:
$$
E=\frac{3}{4} \quad \rightarrow \quad
\AVector{0}{0}{(1,1)}=|0,0>, \quad  \leftarrow \quad \ell_0^{0,(1,1)}=0,
$$
$$
E=\frac{5}{4} \quad \rightarrow \quad
\AVector{0}{0}{(1,2)}=|0,1>, \quad  \leftarrow \quad \ell_0^{0,(1,2)}=0,
$$
$$
E=\frac{7}{4} \quad \rightarrow \quad
\left\{
\begin{array}{l l}
\AVector{1}{-1}{(1,1)}=\frac{1}{\sqrt{2}}\left(|0,2>+i|1,0>\right)
&  \leftarrow \quad \ell_{(-1)}^{1,(1,1)}=-\frac{1}{\sqrt{2}}\\
\AVector{1}{+1}{(1,1)}=\frac{1}{\sqrt{2}}\left(|0,2>-i|1,0>\right)
& \leftarrow \quad \ell_{(+1)}^{1,(1,1)}=\frac{1}{\sqrt{2}},
\end{array}
\right.
$$
$$
E=\frac{9}{4} \quad \rightarrow \quad
\left\{
\begin{array}{l l}
\AVector{1}{-1}{(1,2)}=\frac{1}{\sqrt{2}}\left(|0,3>+i|1,1>\right)
&  \leftarrow \quad \ell_{(-1)}^{1,(1,2)}=-\sqrt{\frac{3}{2}}\\
\AVector{1}{+1}{(1,2)}=\frac{1}{\sqrt{2}}\left(|0,3>-i|1,1>\right)
& \leftarrow \quad \ell_{(+1)}^{1,(1,2)}=\sqrt{\frac{3}{2}},
\end{array}
\right.
$$
$$
E=\frac{11}{4} \rightarrow
\left\{
\begin{array}{l l}
\AVector{2}{-2}{(1,1)}=
\sqrt{\frac{3}{8}}|0,4>+\frac{i}{\sqrt{2}}|1,2>-\frac{1}{\sqrt{8}}|2,0>
&  \leftarrow \quad \ell_{(-2)}^{2,(1,1)}=-2\\
\AVector{2}{0}{(1,1)}=
\frac{1}{2}|0,4>+\frac{\sqrt{3}}{2}|2,0>
&  \leftarrow \quad \ell_{(0)}^{2,(1,1)}=0\\
\AVector{2}{+2}{(1,1)}=
\sqrt{\frac{3}{8}}|0,4>-\frac{i}{\sqrt{2}}|1,2>-\frac{1}{\sqrt{8}}|2,0>
&  \leftarrow \quad \ell_{(+2)}^{2,(1,1)}=+2,
\end{array}
\right.
$$
$$
E=\frac{13}{4} \rightarrow
\left\{
\begin{array}{l l}
\AVector{2}{-2}{(1,2)}=
\frac{\sqrt{5}}{4}|0,5>+\frac{i}{\sqrt{2}}|1,3>-\frac{\sqrt{3}}{4}|2,1>
&  \leftarrow \ell_{(-2)}^{2,(1,2)}=-\sqrt{8}\\
\AVector{2}{0}{(1,2)}=
\sqrt{\frac{3}{8}}|0,5>+\sqrt{\frac{{3}}{8}}|2,1>
&  \leftarrow  \ell_{(0)}^{2,(1,2)}=0\\
\AVector{2}{+2}{(1,2)}=
\frac{\sqrt{5}}{4}|0,5>-\frac{i}{\sqrt{2}}|1,3>-\frac{\sqrt{3}}{4}|2,1>
&  \leftarrow  \ell_{(+2)}^{2,(1,2)}=\sqrt{8}. \\
\end{array}
\right.
$$

The following remarks can now be made:

i) In the basis described by Eqs (15)-(18) it is a trivial matter to
distinguish the states belonging to the same irrep for any $m:n$ ratio,
while in the Cartesian basis this is true only in the 1:1 case.

ii) In the 1:2 case we see that the irreps have degeneracies 1, 1, 2, 2, 3,
3, 4, 4, \dots, i.e. ``two copies'' of the u(2) degeneracies 1, 2, 3, 4,
\dots are obtained.

iii) In the 1:3 cases the degeneracies are 1, 1, 1, 2, 2, 2, 3, 3, 3, \dots,
i.e. ``three copies'' of the u(2) degeneracies are obtained.

iv) It can be easily seen that the 1:n case corresponds to ``n copies'' of
the u(2) degeneracies.

v) Cases with both $m$, $n$ different from unity show more complicated
degeneracy patterns, also correctly reproduced by the above formalism. In
the 2:3 case, for example, the degeneracy pattern is 1, 1, 1, 1, 1, 2, 1, 2,
2, 2, 2, 3, 2, 3, 3, \dots.

vi) The only requirement for each energy eigenvalue to correspond to one
irrep of the algebra is that $m$ and $n$ have to be mutually prime numbers.
If $m$ and $n$ possess a common divisor other than 1, then some energy
eigenvalues will correspond to sums of irreps, i.e. to reducible
representations.

vii) The difference between the formalism used here and the one used in
$^{54,56,57,60}$ is that in the latter case for given $m$ and $n$ appropriate
operators have to be introduced separately for each set of $(p,q)$ values,
while in the present case only one set of operators is introduced.

\subsection{ Connection to W$_3^{(2)}$ }

For the special case $m = 1$, $n=2$ it should be noticed that the  deformed
algebra received here coincides with the finite W algebra  W$_3^{(2)}$
$^{73-76}$. The commutation relations of the W$_3^{(2)}$ algebra are
$$
[H_W, E_W]= 2 E_W, \qquad [H_W, F_W]= -2 F_W, \qquad [E_W,F_W]= H^2_W + C_W,
$$
$$
[C_W,E_W]=[C_W,F_W]=[C_W,H_W]=0,%
$$
while in the $m=1$, $n=2$ case one has the relations
$$
[{\cal N}, {\cal A}^{\dagger}]= {\cal A}^{\dagger}, \qquad [{\cal N}, {\cal A%
}]= -{\cal A}, \qquad [{\cal A},{\cal A}^{\dagger}]= 3 S_0^2 -{\frac{H^2}{4}}
- H S_0 +{\frac{3 }{16}},%
$$
$$
[H, {\cal A}^{\dagger}]= [H, {\cal A}]= [H, S_0] =0,%
$$
with $S_0= {\cal N}+u$ (where $u$ a constant). It is easy to see that the
two sets of commutation relations are equivalent by making the
identifications
$$
F_W= \sigma {\cal A}^{\dagger}, \qquad E_W= \rho {\cal A}, \qquad H_W= -2
S_0 + k H, \qquad C_W= f(H),
$$
with
$$
\rho \sigma = {\frac{4}{3}}, \qquad k={\frac{1}{3}}, \qquad f(H)=-{\frac{4}{9%
}} H^2 +{\frac{1}{4}}.%
$$

\subsection{ Summary}

In conclusion, the two-dimensional anisotropic quantum harmonic oscillator
with rational ratio of frequencies equal to $m/n$, is described dynamically
by a deformed version of the u(2) Lie algebra, the order of this algebra
being $m+n-1$. The representation modules of this algebra can be generated
by using the deformed oscillator algebra. The energy eigenvalues are
calculated by the requirement of the existence of finite dimensional
representation modules. An ``angular momentum'' operator useful for
labelling degenerate states has also been constructed.  The algebras
obtained in the special cases with $1:n$ ratios are shown to correspond to
generalized parafermionic oscillators.
In the special case
of $m:n=1:2$ the resulting algebra has been identified as the finite W
algebra W$_3^{(2)}$.

The extension of the present method to the three-dimensional anisotropic
quantum harmonic oscillator is already receiving attention, since it is of
clear interest in the study of the symmetries underlying the structure of
superdeformed and hyperdeformed nuclei.

%\section{ Acknowledgements}

%One of the authors (DB) has been supported by CEC under contract
%ERBCHBGCT93\-0\-4\-67.

\end{document}